\documentclass{amsart}
\usepackage{amssymb}
\usepackage{epsfig}
\usepackage{amsfonts}
\usepackage{epsfig}

\title{On products of $\mathfrak{sl_n}$ characters and support containment}

\author{Galyna Dobrovolska, Pavlo Pylyavskyy}

\date{July 28th, 2006}

\address{Department of Mathematics, M.I.T., Cambridge, MA 02139}

\email{galyna (at) mit (dot) edu}
\email{pasha (at) math (dot) mit (dot) edu}


\thanks{The research was funded by SPUR program at MIT}

\theoremstyle{plain}
\newtheorem{theorem}{Theorem}

\newtheorem{lemma}[theorem]{Lemma}

\newtheorem{conjecture}[theorem]{Conjecture}
\theoremstyle{definition}
\newtheorem{definition}[theorem]{Definition}
\newtheorem{example}[theorem]{Example}
\theoremstyle{remark}
\newtheorem{remark}[theorem]{Remark}

\def\Imm{\mathrm{Imm}}
\def\TL{\mathit{TL}}

\def\X{\,\,\lower2pt\hbox{
\begin{picture}(0,0)%
\includegraphics{figX.pstex}%
\end{picture}%
\setlength{\unitlength}{1973sp}%
\begingroup\makeatletter\ifx\SetFigFont\undefined%
\gdef\SetFigFont#1#2#3#4#5{%
  \reset@font\fontsize{#1}{#2pt}%
  \fontfamily{#3}\fontseries{#4}\fontshape{#5}%
  \selectfont}%
\fi\endgroup%
\begin{picture}(324,324)(589,-973)
\end{picture}
}}
\def\noXv{\,\,\lower2pt\hbox{
\begin{picture}(0,0)%
\includegraphics{figNoX.pstex}%
\end{picture}%
\setlength{\unitlength}{1973sp}%
\begingroup\makeatletter\ifx\SetFigFont\undefined%
\gdef\SetFigFont#1#2#3#4#5{%
  \reset@font\fontsize{#1}{#2pt}%
  \fontfamily{#3}\fontseries{#4}\fontshape{#5}%
  \selectfont}%
\fi\endgroup%
\begin{picture}(316,316)(293,-969)
\end{picture}
}}

\begin{document}
 
\begin{abstract}
Let $\lambda$, $\mu$, $\nu$ and $\rho$ be dominant weights of $\mathfrak{sl_n}$ satisfying $\lambda + \mu = \nu + \rho$. Let $V_{\lambda}$ denote the highest weight module corresponding to $\lambda$. Lam, Postnikov, Pylyavskyy conjectured a sufficient condition for $V_{\lambda} \otimes V_{\mu}$ to be contained in $V_{\nu} \otimes V_{\rho}$ as $\mathfrak{sl_n}$-modules. In this note we prove a weaker version of the conjecture. Namely we prove that under the conjectured conditions every irreducible $\mathfrak{sl_n}$-module which appears in the decomposition of $V_{\lambda} \otimes V_{\mu}$ does appear in the decomposition of $V_{\nu} \otimes V_{\rho}$.
\end{abstract}

\maketitle

\section{Introduction} \label{intro}

Let $\lambda$, $\mu$ be two dominant weights of $\mathfrak{sl_n}$. Recall that the weight lattice in this case is $\mathbb Z^n/(1, \ldots, 1)$. Thus dominant weights can be viewed as partitions with $n$-th part equal to zero. Equivalently, dominant weights can be associated with Young diagrams with $n-1$ rows. Let $V_{\lambda}$ denote the highest weight module corresponding to $\lambda$. Recall that as $\lambda$ runs through all dominant weights, the $V_{\lambda}$-s constitute the set of irreducible $\mathfrak{sl_n}$-modules. Since $\mathfrak{sl_n}$ is semisimple, the tensor product $V_{\lambda} \otimes V_{\mu} = \bigoplus_{\nu} c_{\lambda, \mu}^{\nu} V_{\nu}$ decomposes into a direct sum of $V_{\nu}$-s. The coefficients $c_{\lambda, \mu}^{\nu}$ which appear in this decomposition are the celebrated Littlewood-Richardson coefficients. Let $\chi_{\lambda}$ denote the polynomial character of irreducible representation $V_{\lambda}$. Then $\chi_{\lambda} = s_{\lambda}(x_1, \ldots, x_n, 0, \ldots)$ is the evaluation of {\it {Schur function}}  $s_{\lambda}$ modulo the relation $x_1 \cdots x_n = 1$. For the background on representation theory of $\mathfrak{sl_n}$ and Schur functions see \cite{Hum}, \cite{Sta}. Schur functions form a basis for the ring $\Lambda$ of symmetric functions with $c_{\lambda, \mu}^{\nu}$ as structure constants. Note that under the substitution $x_{n+1} = \cdots = 0$ the Schur functions $s_{\lambda}$ with $\lambda$ having more than $n$ parts vanish. This causes a subtle difference between multiplication of $\chi_{\lambda}$-s and multiplication of Schur functions: some terms appearing in the latter vanish in the former. 

One can ask when $V_{\lambda} \otimes V_{\mu}$ is contained in $V_{\nu} \otimes V_{\rho}$ as an  $\mathfrak{sl_n}$-module. Of course one way to answer this question is just to say that for all $\kappa$ one should have $c_{\lambda, \mu}^{\kappa} \leq c_{\nu, \rho}^{\kappa}$. However, one might hope to find simple sufficient and/or necessary conditions for the containment to hold. An obviously related question is when the differences of the form $s_{\nu}s_{\rho} - s_{\lambda} s_{\mu}$ are Schur-nonnegative. Some conjectures and results of this form have appeared in the literature, see \cite{BM, FFLP, LLT, LP, LPP, Oko, RS}. 

The first thing to note is that the highest weight appearing in $V_{\lambda} \otimes V_{\mu}$ is $\lambda + \mu$. Thus, in order for $V_{\lambda} \otimes V_{\mu}$ to be a submodule of $V_{\nu} \otimes V_{\rho}$ we need to have $\lambda + \mu \leq \nu + \rho$ in dominance order. It is natural to investigate what happens if we restrict our attention to the case when equality holds, i.e. $\lambda + \mu = \nu + \rho$. For this situation, Lam, Postnikov and Pylyavskyy made a conjecture concerning a sufficient condition for $V_{\lambda} \otimes V_{\mu}$ to be a submodule of $V_{\nu} \otimes V_{\rho}$, or equivalently for $\chi_{\nu}\chi_{\rho} - \chi_{\lambda}\chi_{\mu}$ to be $\chi$-nonnegative. 

Let $\alpha_{ij} = e_i-e_j$ be the roots of the type $A$ root system. Call a polytope {\it {alcoved}} if its faces belong to hyperplanes given by the equations $\langle \alpha_{ij}, \tau \rangle = m$, where $\langle, \rangle$ is the standard inner product and $m \in \mathbb Z$. Alcoved polytopes are studied in \cite{LPo}. Given two weights $\lambda$, $\mu$ one can consider the {\it {minimal alcoved polytope}} $P_{\lambda, \mu}$ containing $\lambda$ and $\mu$. $P_{\lambda, \mu}$ is always a parallelepiped in which $\lambda$ and $\mu$ are a pair of opposite vertices. An example for $\mathfrak{sl_3}$ is shown in Figure \ref{alc3}. The weights $\tau$ inside $P_{\lambda, \mu}$ can be characterized by the following condition: for all $1 \leq i , j \leq n$, the number $\tau_i-\tau_j$ lies weakly between $\lambda_i-\lambda_j$ and $\mu_i-\mu_j$. Let $\nu$ and $\rho$ be another pair of weights.  

\begin{figure} \label{alc3}
  \begin{center}
  \input{minp.pstex_t}
  \end{center}
  \caption{}
\end{figure}

\begin{conjecture} \label{conj}
{\rm \cite{LPP2}} \
If $\lambda + \mu = \nu + \rho$ and $\nu, \rho \in P_{\lambda, \mu}$, then $\chi_{\nu}\chi_{\rho} - \chi_{\lambda}\chi_{\mu}$ is $\chi$-nonnegative.
\end{conjecture}

\begin{example}
It is easy to see in Figure \ref{alc3} that points $\rho = (11,7,0)$, $\nu = (5,2,0)$ lie inside marked $P_{\lambda, \mu}$ with $\lambda = (12,7,0)$, $\mu = (4,2,0)$. In this case $$\chi_{\nu}\chi_{\rho} - \chi_{\lambda}\chi_{\mu} = \chi_{(13,12,0)} + \chi_{(6,4,0)}+\chi_{(7,6,0)}+\chi_{(8,8,0)}+\chi_{(7,3,0)}+\chi_{(8,5,0)}+\chi_{(9,7,0)}$$ $$+\chi_{(10,9,0)}+\chi_{(11,11,0)} + \chi_{(8,2,0)}+\chi_{(9,4,0)}+\chi_{(10,6,0)}+\chi_{(11,8,0)}+\chi_{(12,10,0)}.$$
\end{example}

We prove the following weaker statement.

\begin{theorem} \label{main}
If $\lambda + \mu = \nu + \rho$ and $\nu, \rho \in P_{\lambda, \mu}$, then every $\chi_{\kappa}$ occuring in $\chi_{\lambda} \chi_{\mu}$ with a non-zero coefficient does also occur in $\chi_{\nu}\chi_{\rho}$ with a non-zero coefficient.
\end{theorem}

The paper goes as follows. In Section \ref{hi} we review the theory of Horn-Klyachko inequalities. We prove Lemma \ref{l1} which plays a key role later. In Section \ref{tl} we review Rhoades-Skandera theory of Temperley-Lieb immanants. In Section \ref{proof} we apply the theory of Temperley-Lieb immanants to prove Lemma \ref{l2}. Finally we combine Lemma \ref{l1} and Lemma \ref{l2} to obtain proof of Theorem \ref{main}.

The authors would like to express their gratitude to Thomas Lam and Alex Postnikov, whose work with the second author led to Conjecture \ref{conj}. The authors are also grateful to Denis Chebikin for his superb help with editing the paper.

\section{Horn-Klyachko inequalities} \label{hi}

For a finite set $I=\{i_1>\dots>i_r\}$ of positive integers, define the corresponding partition $\lambda(I)$ by
$$\lambda(I)=(i_1-r, i_2-(r-1),\dots, i_r-1).$$

\begin{definition}
Define $T_{r}^{n}$ to be the set of triples $(I,J,K)$ of subsets of $\{1,\dots,n\}$ of the same cardinality $r$ such that the Littlewood-Richardson coefficient $c_{\lambda(I)\lambda(J)}^{\lambda(K)}$ is positive. A \emph{Horn-Klyachko inequality} for a triple of partitions $\alpha,\beta,\gamma$ has the form 
$$\sum_{k \in K}\gamma_k \leq \sum_{i \in I}\alpha_i+\sum_{j \in J}\beta_j$$
for a triple $(I,J,K)$ in $T_{r}^{n}$ and some $r<n$.
\end{definition}

The following fact was proved in \cite{Kl, KT}, see also \cite{Ful} for a survey:  

\begin{theorem} \label{pos_equiv_horn}
For a triple of partitions $\alpha,\beta,\gamma$ of length $n$, the Littlewood-Richardson coefficient $c_{\alpha\beta}^{\gamma}$ is positive if and only if $\sum_{i=1}^{n}\gamma_i=\sum_{i=1}^{n}\alpha_i+\sum_{i=1}^{n}\beta_i$ and Horn-Klyachko inequalities for $\alpha,\beta,\gamma$ are valid for all $(I,J,K) \in T_{r}^{n}$ and all $r<n$. 
\end{theorem}

Let partitions $\lambda, \mu, \nu, \rho$ with at most $n$ parts satisfy the conditions of Conjecture \ref{conj}, and $\gamma$ be a partiton such that $c_{\lambda\mu}^{\gamma}>0$. Consider a triple $$(I=(i_1,\dots,i_r),J=(j_1,\dots,j_r),K=(k_1,\dots,k_r))$$ in $T_{r}^{n}$. Given permutations $\{l_1,\dots,l_r\}$ of $I$ and $\{m_1,\dots,m_r\}$ of $J$, switch $l_p$ and $m_p$ in some of the pairs  $\{l_p,m_p\}$. This operation yields $2^r$ possible pairs $(I^{\prime}, J^{\prime})$.

\begin{lemma} \label{l1}
Assume there exist permutations $\{l_1,\dots,l_r\}$ of $I$ and $\{m_1,\dots,m_r\}$ of $J$ such that all possible triples $(I^{\prime}, J^{\prime}, K)$ are in $T_{r}^{n}$. Then the Horn-Klyachko inequality corresponding to the triple $(I,J,K)$ holds for $\nu,\rho,\gamma$.
\end{lemma}

\begin{proof}
Since $\nu, \rho \in P_{\lambda, \mu}$, for $i,j \geq 1$ both $\nu_i-\nu_j$ and $\rho_i-\rho_j$ are between $\lambda_i-\lambda_j$ and $\mu_i-\mu_j$, which implies 
$$|(\nu_i-\nu_j)-(\rho_i-\rho_j)| \leq |(\lambda_i-\lambda_j)-(\mu_i-\mu_j)|.$$ 
Rearranging terms, we obtain
$$|(\nu_i+\rho_j)-(\nu_j+\rho_i)| \leq |(\lambda_i+\mu_j)-(\lambda_j+\mu_i)|.$$

This inequality combined with the equality $(\nu_i+\rho_j)+(\nu_j+\rho_i)=(\lambda_i+\mu_j)+(\lambda_j+\mu_i)$ following from $\lambda+\mu=\nu+\rho$, shows that $\nu_i+\rho_j$ and $\nu_j+\rho_i$ are between $\lambda_i+\mu_j$  and $\lambda_j+\mu_i$. We use the fact that for all $i,j \geq 1$ we have 
$$\nu_i+\rho_j \geq \min \{\lambda_i+\mu_j,\lambda_j+\mu_i\}.$$

For every $p$ in $\{1,\dots,r\}$, choose $(l_p^{\prime},m_p^{\prime})$ to be a permutation of $\{l_p,m_p\}$ such that $\lambda_{l_{p}^{\prime}}+\mu_{m_{p}^{\prime}}= \min \{\lambda_{l_p}+\mu_{m_p},\lambda_{m_p}+\mu_{l_p}\}$, and let $I^{\prime}=\{l_1^{\prime},\dots,l_r^{\prime}\}$, $J^{\prime}=\{m_1^{\prime},\dots,m_r^{\prime}\}$ be the corresronding subsets of $\{1,\dots,n\}$. By the assumption of the lemma, $c_{\lambda\mu}^{\gamma}>0$ and $(I^{\prime}, J^{\prime}, K)$ is in $T_{r}^{n}$. Therefore, by Theorem \ref{pos_equiv_horn} the Horn-Klyachko inequality for $\lambda,\mu,\gamma$ and the triple $(I^{\prime},J^{\prime},K)$ holds:
$$\sum_{p=1}^{r}\lambda_{l_{p}^{\prime}}+\sum_{p=1}^{r}\mu_{m_{p}^{\prime}} \geq \sum_{k \in K}\gamma_k.$$  
Observe that 
$$\sum_{i \in I}\nu_i+\sum_{j \in J}\rho_j = \sum_{p=1}^{r}\nu_{l_{p}}+\sum_{p=1}^{r}\rho_{m_{p}} = \sum_{p=1}^{r}(\nu_{l_{p}}+\rho_{m_{p}}) \geq$$
$$\sum_{p=1}^{r}min\{\lambda_{l_p}+\mu_{m_p}, \lambda_{m_p}+\mu_{l_p}\}=\sum_{p=1}^{r}(\lambda_{{l_p}^{\prime}}+\mu_{{m_p}^{\prime}}) \geq \sum_{k \in K}\gamma_k.$$
Therefore, the Horn-Klyachko inequality for $\nu,\rho,\gamma$ and the triple $(I,J,K)$ holds.  \end{proof}

\section{Temperley-Lieb immanants} \label{tl}

In this section we review the theory of Temperley-Lieb immanants developed by Rhoades and Skandera. We limit ourselves to discussing Theorem \ref{th:immdecomp} and Theorem \ref{th:kl} of which we make use in this paper. For detailed exposition of the beautiful results of Rhoades and Skandera we refer reader to the original papers \cite{RS}, \cite{RS2}. One can also find a (more detailed than here) review in \cite{LPP}.

The symmetric functions $h_k = \sum_{i_1\leq \cdots \leq i_k} x_{i_1}\cdots x_{i_k}$ are called the {\it {homogeneous symmetric functions}}. For background on them, see \cite{Sta}. Given two sets $V = (v_1 \geq v_2 \cdots \geq v_n \geq 0)$ and  $U = (u_1 \geq u_2 \cdots \geq u_n \geq 0)$ one can construct the {\it {generalized Jacobi-Trudi matrix}} $X_{V,U} = \left(h_{v_i-u_j}\right)_{i,j = 1}^n$. For example, for $V = (4,3,3,2)$ and $U = (3,2,1,0)$  we get 

$$
X_{V,U}=
\left[
\begin{array}{cccc}
h_1 &  h_2  & h_3 & h_4 \\
1  &  h_1 & h_2 & h_3 \\ 
1 & h_1 & h_2 & h_3 \\
0  &   1  & h_1 & h_2 
\end{array}
\right]
$$

Note that for the operation $\lambda=\lambda(I)=(i_1-r,\ldots,i_r-1)$ defined in Section \ref{hi} we have the Jacobi-Trudi identity $s_{\lambda(I)}=detX_{I,\{r,\ldots,2,1\}}$. (See ~\cite{Sta}).

The {\it {Temperley-Lieb algebra}} $\TL_n(\xi)$ is the $\mathbb
C[\xi]$-algebra generated by $t_1, \ldots, t_{n-1}$
subject to the relations $t_i^2=\xi\, t_i$, $t_i t_j t_i=t_i$ if
$|i-j|=1$, and $t_i t_j = t_j t_i$ if $|i-j| \geq 2$. The dimension of
$\TL_n(\xi)$ equals the $n$-th Catalan number
$C_n=\frac{1}{n+1}\binom{2n}{n}$. 
A {\it 321-avoiding permutation\/} is a permutation $w\in S_n$ that 
has no reduced decomposition of the form $w = \cdots s_i s_j s_i \cdots$ with $|i-j|=1$.
(These permutations are also called {\it fully-commutative.})
A natural basis of the Temperley-Lieb algebra is 
$\{t_w\mid w \textrm{ is a 321-avoiding permutation in } S_n\}$,
where $t_w := t_{i_1} \cdots t_{i_l}$, for a reduced decomposition $w = s_{i_1}
\cdots s_{i_l}$.


For any permutation $v \in S_n$ and a 321-avoiding permutation $w\in S_n$, let $f_{w}(v)$
be the coefficient of the basis element $t_w\in \TL_n(2)$ in the
basis expansion of $(t_{i_1}-1)\cdots
(t_{i_l}-1)\in \TL_n(2)$, where $v=s_{i_1}\cdots s_{i_l}$ is a reduced decomposition.
Rhoades and Skandera~\cite{RS2} defined the {\it Temperley-Lieb
immanant\/} $\Imm_{w}^{\mathrm{TL}}(x)$ of an $n\times n$ matrix $X = (x_{ij})$ by
$$
\Imm_{w}^{\mathrm{TL}}(X):=\sum_{v \in S_n} f_{w}(v)\,x_{1,v(1)} \cdots
x_{n,v(n)}.
$$

\begin{theorem}
\label{th:kl}
{\rm Rhoades-Skandera~\cite[Proposition~2.3, Proposition~3.2]{RS2}} \
Temperley-Lieb immanants of generalized Jacobi-Trudi matrices are Schur-nonnegative.
\end{theorem}

\begin{remark}
In \cite{RS2} two stronger statements (Proposition~2.3 and Proposition~3.2) are proved, from which Theorem \ref{th:kl} follows in a straightforward way.
\end{remark}

A product of generators (decomposition) 
$t_{i_1}\cdots t_{i_l}$ in the Temperley-Lieb algebra $\TL_n$ can be
graphically presented by a {\it Temperley-Lieb diagram\/} with $n$ non-crossing
strands connecting the vertices $1,\dots,2n$, possibly with
some internal loops. The left endpoints are assumed to
be labeled $1,\dots, n$ from top to bottom and the right endpoints 
are assumed to be labeled $2n,\dots, n+1$ from top to bottom.

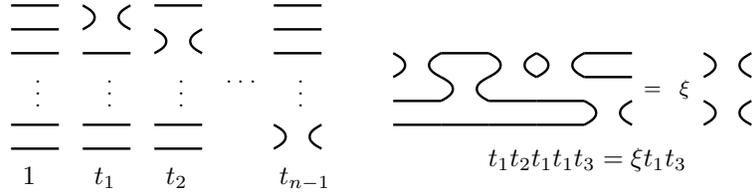
\begin{figure}[h!]
     \begin{center}
     \input{figvam9.pstex_t} 
     \end{center}
\end{figure}

The map that sends $t_w$ to the non-crossing matching
given by its Temperley-Lieb diagram is a bijection between
basis elements $t_w$ of $\TL_n$,  where $w$ is 321-avoiding,
and non-crossing matchings on the vertex set $[2n]$. 

Following~\cite{RS2}, for a subset $S\subset [2n]$, 
let us say that a Temperley-Lieb
diagram (or the associated element in $\TL_n$) is {\it $S$-compatible\/} 
if each strand of the diagram has one endpoint in $S$ and the
other endpoint in its complement $[2n]\setminus S$.  
Coloring vertices in $S$ black 
and the remaining vertices
white, a basis element $t_w$ is $S$-compatible if and only if each edge 
in the associated matching has two vertices of different color.  
Let $\Theta(S)$ denote the set of all 321-avoiding permutations
$w\in S_n$ such that $t_w$ is $S$-compatible. An example for $n=5$, $S = \{3,6,7,8,10\}$ is shown in the figure below, where all possible compatible non-crossing matchings are presented.

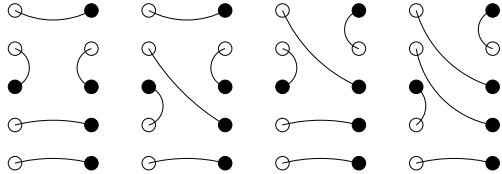
\begin{figure}[h!]
    \begin{center}
    \input{figvam10.pstex_t}
    \end{center}
\end{figure}

For two subsets $I, J\subset [n]$ of the same cardinality, let $\Delta_{I,J}(X)$
denote the {\it minor\/} of an $n\times n$ matrix $X$ in the row set $I$ and the
column set $J$.
Let $\bar I := [n]\setminus I$ and let $I^\wedge := \{2n+1-i\mid i\in I\}$.

\begin{theorem} 
\label{th:immdecomp}
{\rm Rhoades-Skandera~\cite[Proposition~4.4]{RS2}, cf.~Skandera~\cite{Ska}} \
For two subsets $I,J\subset [n]$ of the same cardinality
and $S=J\cup (\bar I)^\wedge$, we have
$$
\Delta_{I,J}(X)\cdot \Delta_{\bar I, \bar J}(X)
= \sum_{w \in \Theta(S)} \Imm_w^{\mathrm{TL}}(X).
$$
\end{theorem}

\begin{example} 
Take $I = \{1,2\}$, $J = \{1,3\}$, and 
$$
X = 
\left[
\begin{array}{cccc}
x_{11} &  x_{12}  & x_{13} & x_{14}\\
x_{21}  &  x_{22} & x_{23} & x_{24}\\
x_{31} & x_{32} & x_{33} & x_{34}\\
x_{41}  &   x_{42}  & x_{43} & x_{44}
\end{array}
\right]
$$
Then $S = \{1,3,5,6\}$, and the elements of $\Theta(S)$ are shown in the figure below. 

\begin{figure}[h!]
     \begin{center}
     \input{figvam13.pstex_t}
     \end{center}
\end{figure}

In this case Theorem \ref{th:immdecomp} yields the decomposition

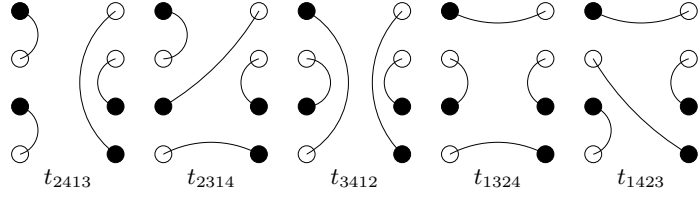
\begin{eqnarray*}
\left|
\begin{array}{cc}
x_{11} &  x_{13}\\
x_{21} & x_{23}
\end{array}
\right|
\times
\left|
\begin{array}{cc}
x_{32} &  x_{34}\\
x_{42} & x_{44}
\end{array}
\right|
= \sum_{w \in \Theta(S)} \Imm_w^{\mathrm{TL}}(X).
\end{eqnarray*}
\end{example}

\section{Proof of the main theorem} \label{proof}

We begin by proving the following lemma, which shows that the condition of Lemma \ref{l1} can be achieved. 

\begin{lemma} \label{l2}
In the setup of Lemma \ref{l1}, there exist permutations $\{l_1,\dots,l_r\}$ and $\{m_1,\dots,m_r\}$ of $I$ and $J$ respectively such that all possible triples $(I^{\prime}, J^{\prime}, K)$ are in $T_{r}^{n}$.
\end{lemma}

\begin{proof}
Let $X_{V,U}$ be the generalized Jacobi-Trudi matrix for column set $U=(r,r,r-1,r-1,\ldots,1,1)$, and row set $V=I \cup J$ in some chosen non-increasing arrangement. Let $\#I$ and $\#J$ denote the sets of numbers of the rows of $I$ and $J$ in the chosen non-increasing arrangement of $I \cup J$. Since $(I,J,K) \in T_{r}^{n}$, we have $c_{\lambda(I)\lambda(J)}^{\lambda(K)}>0$. Hence $s_{\lambda(K)}$ is present in the decomposition of  $s_{\lambda(I)} s_{\lambda(J)}$, which by Jacobi-Trudi identity equals to the product $\Delta_{\#I,\{2r,2r-2,\ldots,2\}} \Delta_{\#J,\{2r-1,2r-3,\ldots,1\}}$ of complementary minors of $X_{V,U}$. This product, in turns, by Theorem \ref{th:immdecomp} equals to $\sum_{w \in \Theta(S)} \Imm_w^{\mathrm{TL}}(X_{V,U})$, where $S=\#J \cup \{4r,4r-2,\ldots,2r+2\}$ is the subset of the vertices ${1,\ldots,4r}$ of the Temperley-Lieb diagram which are colored black.

Since $s_{\lambda(K)}$ is in the Schur function decomposition of $\sum_{w \in \Theta(S)} \Imm_w^{\mathrm{TL}}(X_{U,V})$, it is present in the Schur function decomposition of one of the immanants $\Imm_w^{\mathrm{TL}}(X_{V,U})$ for some $321$-avoiding permutation $w \in \Theta(S)$. For this $321$-avoiding permutation $w$, the basis element $t_w$ and the corresponding non-crossing matching $M_w$ of the Temperley-Lieb diagram with columns $V$ and $U$ are $S$-compatible. Therefore, all edges of $M_w$ have endpoints of different color in the Temperley-Lieb diagram on vertices $\{1,2,\ldots,4r\}$ where $S$ is colored black and $[4r]/S$ colored white.  

We proceed now to construct the needed permutations $\{l_1,\ldots,l_r\}$ of $I$ and $\{m_1,\dots,m_r\}$ of $J$ based on $S$ and $M_w$. We go along $V$ from top to bottom (see Figure \ref{horn_sec4_fig1}$(i)$) and label vertices in $I$ that are connected to vertices in $J$ by edges in $M_w$ (suppose that there are $k$ such vertices in $I$) with variables $l_1,\ldots,l_k$ as we meet them. We also label the vertex in $J$ connected to $l_i$ $(i \leq k)$ by $m_i$.

Next, we remove the vertices $l_1,\ldots,l_k,m_1,\ldots,m_k$ from $V$ and call the remaining set $V^-$. We also go along $U$ and discard every pair of vertices in $U$ connected by an edge in $M_w$, and call the remaining set $U^-$. We go along $V^-$ from top to bottom and label the white vertices that we meet by $l_{k+1},\ldots,l_r$, and the black vertices we meet by $m_{k+1},\ldots,m_r$ from top to bottom. For $f \geq 1$, we also label the vertices in $U^-$ connected by edges in $M_w$ to $l_{k+f}$ by $p_{k+f}$, and those connected to $m_{k+f}$ by $q_{k+f}$. (See Figure \ref{horn_sec4_fig1}$(ii)$). Note that every vertex in $V$ between adjacent vertices of $V^{\prime}$ is connected by an edge in $M_w$ to another vertex between the same vertices of $V^-$ because $M_w$ is a non-crossing, and the same is true about $U$. Therefore, in building $V^-$ and $U^-$ we discarded segments of even lengths from $V$ and $U$.  

\emph{Claim.} For $f \geq 1$, vertices $l_{k+f}$ and $q_{k+f}$ are white and odd-numbered in the Temperley-Lieb diagram for $S$ and $M_w$; vertices $p_{k+f}$ and $m_{k+f}$ are black and even-numbered. Also, $l_{k+f+1}>m_{k+f}>l_{k+f}$ and $p_{k+f+1}<q_{k+f}<p_{k+f}$. (See Figure \ref{horn_sec4_fig1}$(ii)$)

\emph{Proof.} Since we discarded segments of even lengths from $U$ to obtain $U^-$ and the colors in $U$ were alternating from top to bottom beginning with the black even vertex $4r$, the colors in $U^-$ are also alternating from top to bottom beginning with a black even vertex. Therefore, vertices in $U^-$ from top to bottom are $p_{k+1}>q_{k+1}>p_{k+2}>\ldots>p_{r}>q_{r}$, where $p_{k+f}$ is black and $q_{k+f}$ is white for $f \geq 1$. Because the restriction of the matching $M_w$ to $U^- \cup V^-$ is non-crossing, the inequalities $p_{k+1}>q_{k+1}>p_{k+2}>\ldots>p_{r}>q_{r}$ for $U^-$ imply that $l_{k+1}<m_{k+1}<l_{k+2}<\ldots<l_{r}<m_{r}$ for $V^-$. The colors in $V^-$ alternate and have a white odd vertex at the top because the colors in $U^-$ alternate with a black even vertex at the top. Therefore, $l_{k+f}$ is white and $m_{k+f}$ is black for $f \geq 0$. The statements about being odd/even now follow from the fact that we discarded segments of even lengths from $U$ and $V$ to obtain $U^-$ and $V^-$. 

We now build a new coloring $S^{\prime}$ of $U \cup V$ based on the transpositions $(l_p,m_p)$ that may have occured in going from $I,J$ to $I^{\prime},J^{\prime}$. We only allow ourselves to recolor both elements in a pair $\{2m,2m-1\} \in \{2r+1,\ldots,4r\}$ of vertices in the second column of Temperley-Lieb diagram for $S_0^{\prime}=\#J^{\prime} \cup \{4r,4r-2,\ldots,2r+2\}$, because the columns $4r+1-2m$ and $4r+2-2m$ of $X_{U,V}$ are identical and hence such a recoloring produces the same pair of complementary minors $\Delta_{\#I^{\prime},\{4r-1,\ldots,2r+1\}}$, $\Delta_{\#J^{\prime},\{4r,\ldots,2r+2\}}$ of $X_{V,U}$ as $S_0^{\prime}$ does, and therefore by Jacobi-Trudi identity the product of these complementary minors is $s_{\lambda(I^{\prime})}s_{\lambda(J^{\prime})}$.

\emph{Rule of recoloring.} For every pair $l_{k+f}$ and $m_{k+f}$ $(f \geq 1)$ that exchanged colors in transition from $I,J$ to $I^{\prime}, J^{\prime}$, recolor the pairs $(p_{k+f}, p_{k+f}-1),(p_{k+f}-2,p_{k+f}-3),\ldots,(q_{k+f}+1,q_{k+f})$. The recoloring is permissible because the vertex $p_{k+f}$ is even by the Claim. (See Figure \ref{horn_sec4_fig1}$(iii)$)

\begin{figure}[hbtp]
    \begin{center}
    \input{horn_sec4_fig1.pstex_t}
    \end{center}
    \caption{}  \label{horn_sec4_fig1}
\end{figure}

\emph{Why the rule produces a coloring compatible with $M_w$.} The vertices between $p_{k+f}$ and $q_{k+f}$ either all changed color or all stayed the same, so an edge in $M_w$ that connected two vertices in $U$ between $p_{k+f}$ and $q_{k+f}$ now has its endpoints changed or not changed simultaneously, so they are of different color in the new coloring.

A pair $(l_{k+f},m_{k+f})$ changes color simultaneously with the pair $(p_{k+f},q_{k+f})$, so $l_{k+f}$ and $p_{k+f}$ , and $m_{k+f}$ and $q_{k+f}$ change or do not change their color simultaneously, so the endpoints of the edges between $U^-$ and $V^-$ remain colored differently in the new coloring. 

A pair of vertices $(l_p,m_p)$ in $V$ connected by an edge in $M_w$ changes color simultaneously when the corresponding transposition occurs, so the endpoints of such an edge remain colored differently. Finally, a pair of vertices in $U$ between $q_{k+f}$ and $p_{k+f+1}$ connected by an edge in $M_w$ never changes color, so such an edge has its endpoints colored differently in the new coloring. We considered all possibilities for an edge in $M_w$ relative to $U^-$ and $V^-$ in a non-crossing matching, so $M_w$ is compatible with the new coloring.   

We already noticed that the new coloring produces the product of complementary minors of $X_{V,U}$ equal to $s_{\lambda(I^{\prime})} s_{\lambda(J^{\prime})}$. The fact that the new coloring is compatible with $M_w$ implies that the immanant $\Imm_w^{\mathrm{TL}}(X_{V,U})$ is present in the decomposition $s_{\lambda(I^{\prime})} s_{\lambda(J^{\prime})}=\sum_{w \in \Theta(S^{\prime})} \Imm_w^{\mathrm{TL}}(X_{V,U})$. Since $s_{\lambda(K)}$ is in the decomposition of $\Imm_w^{\mathrm{TL}}(X_{V,U})$ which is Schur-nonnegative by Theorem \ref{th:kl}, $s_{\lambda(K)}$ is present in the Schur function decomposition of $s_{\lambda(I^{\prime})} s_{\lambda(J^{\prime})}$. Therefore $c_{\lambda(I^{\prime})\lambda(J^{\prime})}^{\lambda(K)}>0$ and $(I^{\prime},J^{\prime},K) \in {T_r}^n$ for all $I^{\prime},J^{\prime}$ that can be obtained by transposing pairs $(l_p,m_p)$ in $I$,$J$.
\end{proof}

We are ready to prove Theorem \ref{main}.

\begin{proof}
From Lemma \ref{l1} and Lemma \ref{l2} it follows that whenever the Horn-Klyachko inequality for triple $(I,J,K)$ holds for $\lambda, \mu, \gamma$, it also holds for $\nu,\rho,\gamma$. Thus all possible $\gamma$-s for which all needed Horn-Klyachko inequlities hold or, equivalently, $c_{\lambda, \mu}^{\gamma} > 0$, also have the property that $c_{\nu, \rho}^{\gamma} > 0$. 
\end{proof}

\end{document}